\documentclass[10pt,a4paper,english]{article} 
\usepackage{color}
\usepackage[utf8]{inputenc}

\usepackage{amsmath, amssymb, latexsym}
\usepackage[english]{babel}
\usepackage{float}
\usepackage{setspace}
\usepackage{stmaryrd}
\usepackage{enumerate}
\usepackage{eufrak}
\usepackage{graphicx}

\def\CA{{\mathcal A}}

\def\CD{{\mathcal D}}

\def\CH{{\mathcal H}}

\def\CP{{\mathcal P}}

\def\pa{{\partial}}

\newcommand{\sca}[2]{\langle#1,#2\rangle}
\newcommand{\barr}{\overline}

\newcommand{\beq}{\begin{equation}}
\newcommand{\eeq}{\end{equation}}

\newcommand{\Supp}{\textrm{Supp~}}

\renewcommand{\Re}{{\rm Re\,}}
\renewcommand{\Im}{{\rm Im\,}}
\newtheorem{theorem}{Theorem}[section]
\newtheorem{lemma}[theorem]{Lemma}

\newtheorem{proposition}[theorem]{Proposition}

\newtheorem{corollary}[theorem]{Corollary}

\makeatletter
\@addtoreset{equation}{section}

\makeatother
\title{Spectral projections of the complex cubic oscillator}
\author{Raphaël \textsc{Henry}\footnote{
D\'epartement de Math\'ematiques, Batiment 425, Universit\'e Paris Sud, 91405 Orsay Cedex, France.
email: raphael.henry@math.u-psud.fr}
\thanks{The author is supported  by the ANR NOSEVOL.}}
\date{}
\begin{document}
\maketitle
\begin{abstract}
 We prove the spectral instability of the complex cubic oscillator $-\frac{d^2}{dx^2}+ix^3+i\alpha x$ for non-negative values
of the parameter $\alpha$, by getting the exponential growth rate of $\|\Pi_n(\alpha)\|$, where $\Pi_n(\alpha)$ is the spectral
projection associated with the $n$-th eigenvalue of the operator. More precisely, we show that for all non-negative $\alpha$
\[
 \lim\limits_{n\to+\infty}\frac{1}{n}\log\|\Pi_n(\alpha)\| = \frac{\pi}{\sqrt{3}}.
\]
\end{abstract}

\section{Introduction}
\begin{figure}[H]\begin{center}
 \includegraphics[scale=0.4]{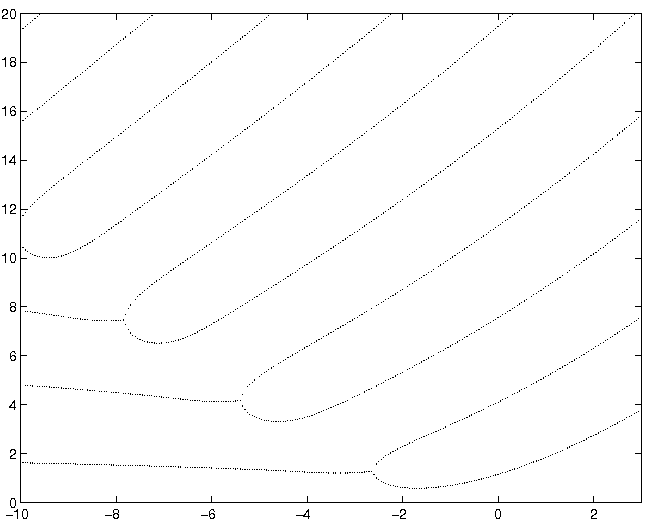}
 \caption{Real parts of the eigenvalues of $\CA_\alpha$ as functions of $\alpha$. Each pair of
consecutive eigenvalues becomes non-real, complex conjugate on the left of the branch point.}\label{figSpCub}
\end{center}
\end{figure}

We consider the complex cubic oscillator

\beq\label{defAalpha}
 \CA_\alpha = -\frac{d^2}{dx^2}+ix^3+i\alpha x,~~\alpha\in\mathbb{R}
\eeq
on the real line. We define $\CA_\alpha$ by extension of the operator
\[
  \CA_\alpha^0 = -\frac{d^2}{dx^2}+ix^3+i\alpha x,~~~~~
\CD(\CA_\alpha^0) = \mathcal{C}_0^\infty(\mathbb{R}),
 \]
which is accretive, so we can define $\CA_\alpha := \barr{\CA_\alpha^0}$ as its closure.
$\CA_\alpha$ is then maximally accretive, with domain
\[
 \CD(\CA_\alpha) = H^2(\mathbb{R})\cap L^2(\mathbb{R} ; x^6dx).
\]
The cubic oscillator presented here has been studied in \cite{DelTri} and \cite{Trinh1}. 
It also belongs to the class of operators considered in
\cite{Shin}. Let us mention \cite{GreMaiMar} as well, which deals with a quadratic perturbation of the cubic $ix^3$ potential.\\
The operator $\CA_\alpha$ has compact resolvent, and its eigenvalues 
$(\lambda_n(\alpha))_{n\geq1}$
are simple in the sense of the geometric multiplicity.\\
The properties of the complex cubic oscillator and its variants (the potential $x^2+ix^3$, for instance),
have been widely studied in the
past few years (see \cite{Ben, BenBoe, CalGraMai, CalMai, DelPhaI, DelPhaII, DelTri, GreMaiMar, KreSie, Trinh1, Trinh2, Shin}).
As a non-selfadjoint operator, it has a surprising property: its spectrum is purely real for $\alpha\geq0$
(see \cite{BenBoe} for numerical observations and \cite{Shin} for a rigorous proof).
This property is suspected to be related with the so-called \emph{$\mathcal{PT}$-symmetry} of the operator,
namely
\[
\mathcal{PT}\CA_\alpha = \CA_\alpha\mathcal{PT},
\]
where $\mathcal{P}$ and $\mathcal{T}$, denoting respectively the spatial symmetry and time inversion operators, act as follows:
\[
 (\CP u)(x) = u(-x)~~~\textrm{ and }~~~(\mathcal{T}u)(x) = \barr{u(x)}.
\]
The complex cubic oscillator is a toy model in the study of $\mathcal{PT}$-symmetric operators.\\
One of the main questions arising from this property of real spectrum is the following: does $\CA_\alpha$ share some 
other similarities with selfadjoint operators? More precisely, does the family of eigenfunctions form a basis of $L^2(\mathbb{R})$ in some sense?
Is the spectrum stable under perturbations of the operator? What can one say about the behavior of the eigenvalues
for negative values of $\alpha$?
Some of these questions have already been answered, while other have been stated as conjectures.
For instance, it has been established in \cite{KreSie} that the eigenfunctions of $\CA_\alpha$ do not form a Riesz basis, as well as the existence of 
non-trivial pseudospectra.\\
The properties of the spectrum of $\CA_\alpha$ for negative $\alpha$ have not been completely understood yet.
Numerical simulations
(see \cite{DelPhaI}, \cite{DelPhaII}, \cite{DelTri}), reproduced on Figure
\ref{figSpCub}, suggest that, for any $n\geq1$, there exists a critical value 
$\alpha_n^{crit}<0$ of the parameter such that $\lambda_n(\alpha)$ is real for $\alpha>\alpha_n^{crit}$. For
$\alpha = \alpha_n^{crit}$, $\lambda_n(\alpha_n^{crit})$ seems to cross an adjacent eigenvalue, forming
for $\alpha<\alpha_n^{crit}$ a complex conjugate pair lying away from the real axis.
Regarding the analysis for large eigenvalues which we will perform in the following,
the simulation suggests that, for any fixed $\alpha<0$, the eigenvalues $\lambda_n(\alpha)$ are real for $n$ large enough,
but it does not seem to be proved yet. Therefore, we will only consider non-negative values of $\alpha$ in the following.\\

Our goal is to measure the spectral instability of the operator $\CA_\alpha$.
As mentioned above, the instability of the eigenvalues $\lambda_n(\alpha)$ has already been highlighted in 
\cite{KreSie} by proving the existence of non-trivial pseudospectra. We now want to understand more accurately this phenomenon,
following the approach of \cite{Dav2}, \cite{DavKui} and \cite{Hen}.\\
To this purpose, we define the instability indices
\beq\label{defKappaCub}
\kappa_n(\alpha) = \|\Pi_n(\alpha)\|,
\eeq
where $\Pi_n(\alpha)$ denotes the spectral projection of $\CA_\alpha$ associated with the eigenvalue $\lambda_n(\alpha)$ (the eigenvalues
being labelled in increasing order).
We shall first consider the question of algebraic multiplicity for the eigenvalues $\lambda_n(\alpha)$,
that is, whether there exist associated Jordan blocks or not. The algebraic simplicity of the eigenvalues has
been proved for all $n\geq1$ in \cite{GreMaiMar} in the case of a potential of the form
$ax^2+i\sqrt{\beta}x^3$.  Here, by an independent proof, we shall get the algebraic simplicity of $\lambda_n(\alpha)$, but only for $n$ large enough,
which will be enough to achieve the proof of our main statement.
Hence, for $n$ large enough, the expression
\beq\label{kappan}
 \kappa_n(\alpha) = \frac{\|u_n^\alpha\|^2}{|\sca{u_n^\alpha}{\bar u_n^\alpha}|}
\eeq
will hold,
where $u_n^\alpha$ denotes an eigenfunction of $\CA_\alpha$ associated with the eigenvalue $\lambda_n(\alpha)$
(see \cite{AslDav}). We will use this formula to prove the following theorem, which is the main statement of our work.
\begin{theorem}\label{mainCubic}
For all $\alpha\geq0$, we have
\beq
\lim\limits_{n\to+\infty}\frac{1}{n}\log\kappa_n(\alpha) = \frac{\pi}{\sqrt{3}}.
\eeq
\end{theorem}
Let us recall that the same question was considered in \cite{Dav2,DavKui,Hen} in the case of anharmonic oscillators
$-\frac{d^2}{dx^2}+e^{i\theta}|x|^m$, $m>0$, $|\theta|<\min\{(m+2)\pi/4,(m+2)\pi/2m\}$. More precisely, it has been proved that the spectral projections of
these operators grow faster than any power of $n$ as $n\rightarrow\infty$ \cite{Dav2}, and the exponential growth rate was precisely obtained
for $m=2$ in \cite{DavKui} and for every even exponent $m$ in \cite{Hen}.\\
The proof of Theorem \ref{mainCubic} lies on WKB estimates of the eigenfunctions in the complex plane. This method has already been used
in \cite{Hen} in the even anharmonic case. However, here we will have to manage the sub-principal
term $i\alpha x$ in the potential.\\
Some results from \cite{KreSie} can be recovered immediately from Theorem \ref{mainCubic}:
\begin{corollary}
 For all $\alpha\geq0$, the eigenfunctions of $\CA_\alpha$ do not form a Riesz basis.
\end{corollary}
\textbf{Proof: } 
Let $(u_n^\alpha)_{n\geq1}$ be a family of eigenfunctions for $\CA_\alpha$ associated with the eigenvalues $(\lambda_n(\alpha))_{n\geq1}$.
Let us recall that $(u_n^\alpha)_{n\geq1}$ is said to be a Riesz basis if it spans a dense subset of $L^2(\mathbb{R})$ and if there exists
$C>0$ such that, for all $\phi\in L^2(\mathbb{R})$,
\beq\label{inegRieszDef}
 C^{-1}\sum_{n=1}^{+\infty}|\sca{\phi}{u_n^\alpha}|^2\leq\|\phi\|^2\leq C\sum_{n=1}^{+\infty}|\sca{\phi}{u_n^\alpha}|^2\,.
\eeq
According to Lemma \ref{propII} and Proposition \ref{propDenomCubic} (which provides algebraic simplicity for large eigenvalues of $\CA_\alpha$),
we can choose the eigenfunctions $u_n^\alpha$ such that, for $n,m\geq1\,$ and $n$ large enough, 
$\sca{u_n^\alpha}{\barr{u_m^\alpha}} = \delta_{n,m}$. Hence according to (\ref{exprII}), we have
$\kappa_n(\alpha) = \|u_n^\alpha\|^2$ for $n$ large enough.
Using that $\kappa_n(\alpha)\rightarrow+\infty$ as $n\rightarrow+\infty$,
it is then straightforward to check that the sequence $\phi_n = \barr{u_n^\alpha}$ can not satisfy (\ref{inegRieszDef}). 
\hfill $\square$\\
Furthermore, the pseudospectra in the neighborhood of an eigenvalue are known to grow proportionally to the corresponding instability index
(see \cite{AslDav}, \cite{TrEm}). Hence the exponential growth obtained in Theorem \ref{mainCubic}
enables us to confirm the presence of nontrivial pseudospectra \cite{KreSie}, and to somehow describe its shape near the eigenvalues. 
\\

Section \ref{sAsymp} is devoted to the estimates on the eigenfunctions needed to prove Theorem \ref{mainCubic}. 
The proof itself is achieved in Section \ref{sCcl}.

\section{Asymptotic behavior of the eigenfunctions}\label{sAsymp}
\subsection{Preliminary scale change}\label{pScaling}
Let us first perform the following scale change.
Let us recall that for all $\alpha\geq0$,
the spectrum of $\CA_\alpha$ is real, and let us denote the eigenvalues, labelled in increasing order, by $\lambda_n(\alpha)$. We set
\beq\label{defHnx}
\left\{
 \begin{array}{ccc}
  h_n & = & \lambda_n(\alpha)^{-5/6} \\
\tilde x & = & h_n^{2/5}x.
 \end{array}
\right.
\eeq
The operator $(\CA_\alpha-\lambda_n(\alpha))$ then writes
\[
 -h_n^{4/5}\frac{d^2}{d\tilde x^2}+ih_n^{-6/5}\tilde x^3+i\alpha h_n^{-2/5}\tilde x-\lambda_n(\alpha) =
 h_n^{-6/5}\left(-h_n^2\frac{d^2}{d\tilde x^2}+i\tilde x^3+i\alpha h_n^{4/5}\tilde x-1\right),
\]
and we are reduced to the study of the kernel of
\[
 \CA_\alpha(h) = -h^2\frac{d^2}{dx^2}+ix^3+i\alpha h^{4/5}x-1.
\]
An eigenfunction $u_n^\alpha$ of $\CA_\alpha$ associated with $\lambda_n(\alpha)$ can be written as
\beq\label{UnUh}
 u_n^\alpha(x) = \psi_\alpha(h_n^{2/5}x,h_n)=\psi_\alpha(\lambda_n(\alpha)^{-1/3}x,h_n),
\eeq
where $\psi_\alpha(\cdot,h_n)$ is a solution of 
\beq\label{equaPsi}
 \CA_\alpha(h_n)\psi_\alpha(\cdot,h_n) = 0,~~~\psi_\alpha(\cdot,h_n)\in L^2(\mathbb{R}).
\eeq
Notice that the condition $\psi_\alpha(\cdot,h_n)\in L^2(\mathbb{R})$, together with (\ref{equaPsi}),
ensure that $\psi_\alpha(\cdot,h_n)$ belongs to the domain
$\CD(\CA_\alpha(h_n)) = \CD(\CA_\alpha)$ (see for instance Theorem \ref{CorWKB} below).
Thus, we will now work on these solutions $\psi_\alpha$.\\
From now on, $\alpha$ is assumed to be fixed and non-negative.

\subsection{Behavior of the eigenfunctions away from the turning points}\label{pWKB}
\begin{figure}[H]\begin{center}
 \scalebox{0.3}{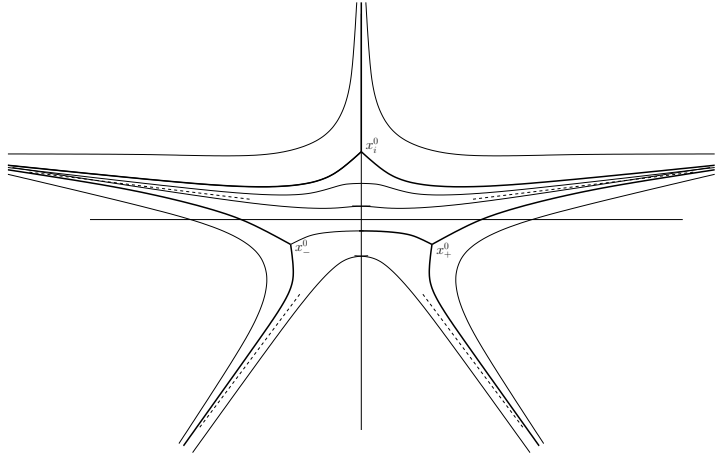}
 \caption{Stokes lines of the operator $\CA_0 = -h^2\frac{d^2}{dx^2}+ix^3-1$. 
Bold lines are those starting from the turning points. Dashed lines are the asymptotic directions $D_k$, $k=0,\dots,4$.}
\label{figStokes}
\end{center}
\end{figure}

In this subsection, we determine the global asymptotic behavior of the solutions $\psi_\alpha(x,h)$ of
\beq\label{probCubic}
 \CA_\alpha(h)\psi_\alpha(x,h) = 0,~~~\psi_\alpha(\cdot,h)\in L^2(\mathbb{R})
\eeq
as $h\rightarrow0$.\\
More precisely, we want to understand the behavior of $\psi_\alpha$ in a domain of the complex plane avoiding the zeroes
(called turning points of the equation) of the potential
\[
 V_\alpha(x,h) = ix^3+i\alpha h^{4/5}x-1.
\]
Let $x_+^\alpha(h)$, $x_-^\alpha(h)$ and $x_{\mathbf{i}}^\alpha(h)$ denote the zeroes of $V_\alpha(\cdot,h)$,
respectively starting at $h=0$ from the zeroes
$x_+^0 = e^{-i\pi/6}$, $x_-^0 = e^{-5i\pi/6}$ and $x_{\mathbf{i}}^0 = i$ of the potential
\[
 V_0(x) = ix^3-1.
\]
Note that for $h$ small enough, $x_\pm^\alpha(h)$, $x_{\mathbf{i}}^\alpha(h)$ are simple zeroes of $V_\alpha(\cdot,h)$.\\
To understand the asymptotic properties of the solutions of (\ref{probCubic}), 
it will be useful to analyze the geometry of the level curves (\emph{Stokes lines}) of the function
\[
 x\mapsto\Re\int_{x_+^\alpha(h)}^x\sqrt{V_\alpha(z,h)}~dz,
\]
where $\sqrt{V_\alpha}$ is holomorphic in
\[
 D_h^\alpha = \mathbb{C}\setminus\bigcup_{\tiny{\sigma\in\{+,-,\mathbf{i}\}}}\{(1+r)x_\sigma^\alpha(h) : r>0\},
\]
and $\sqrt{V_\alpha(0,h)} = i$.\\The path of integration is included in $D_h^\alpha$.\\
Let us notice that $x_+^\alpha(h)$ and $x_-^\alpha(h)$ belong to a common, bounded Stokes line, joining the two points:
\[
 \Re\int_{x_-^\alpha(h)}^{x_+^\alpha(h)}\sqrt{V_\alpha(z,h)}~dz = 0.
\]
Let us denote this line by $\ell_f^\alpha(h)$. It is the only bounded Stokes line for $\CA_\alpha$
(see Figure \ref{figStokes}).\\
On the other hand, there are seven unbounded Stokes lines starting from $x_\pm^\alpha(h)$, $x_{\mathbf{i}}^\alpha(h)$, with the 
five asymptotic directions as $|x|\rightarrow+\infty$, 
\[
D_k = \arg^{-1}\left\{\frac{\pi}{10}+\frac{2k\pi}{5}\right\},~~~k=0,\dots,4.
\]
Among those Stokes lines, one 
is starting from $x_{\mathbf{i}}^\alpha(h)$ and has asymptotic direction
$D_1 = i\mathbb{R}^+$ ; let us denote it by $\ell_{\mathbf{i}}^\alpha(h)$. Notice that for $h=0$, $\ell_{\mathbf{i}}^\alpha(0) = i[1,+\infty[$.\\
For $\varepsilon>0$, let
\beq\label{NotLdS}
\ell_{f,\varepsilon}^0 = \{x\in\mathbb{C} : d(x,\ell_f^0(0))<\varepsilon\},
\eeq
and
\beq\label{NotLdSf}
 \ell_{i,\varepsilon}^0 = \{x\in\mathbb{C} : d(x,i[1,+\infty[)<\varepsilon\}.
\eeq
Hence, for all $\varepsilon>0$ fixed, there exists $h_0>0$ such that, for all $h\in]0,h_0[$,
\beq\label{incluLdS}
\ell_f^\alpha(h)\subset\ell_{f,\varepsilon}^0,~~~\ell_{\mathbf{i}}^\alpha(h)\subset\ell_{i,\varepsilon}^0.
\eeq
Finally, let
\beq
\Gamma_\varepsilon = \mathbb{C}\setminus(\ell_{f,\varepsilon}^0\cup\ell_{i,\varepsilon}^0).
\eeq
In the following theorem, $(h_n)_{n\geq1}$ is the sequence defined in (\ref{defHnx}).
\begin{theorem}\label{CorWKB}
 Let $\varepsilon>0$ be fixed. There exists $N\geq1$ such that, for all $n\geq N$, there exists a unique solution
$\psi_1^\alpha(x,h_n)\in L^2(\mathbb{R})$ of
\beq\label{solPsin}
 \CA_\alpha(h_n)\psi_1^\alpha(\cdot,h_n) = 0
\eeq
satisfying
\beq\label{WKBcubXn}
\psi_1^\alpha(x,h_n) = \frac{e^{-i\pi/8}}{x^{3/4}}(1+o(1))
\exp\left(-\frac{1}{h_n}\int_{x_+^\alpha(h_n)}^x\sqrt{V_\alpha(z,h_n)}~dz\right)
\eeq
as $|x|\rightarrow+\infty$ in $\Gamma_\varepsilon$, uniformly with respect to $n\geq N$.\\
Moreover, there exists a sequence $(u_j^\alpha)_{j\geq1}$ of functions, holomorphic on $\Gamma_\varepsilon$,
such that, for every $j_0\geq1$ and $x\in\Gamma_\varepsilon$,
\begin{eqnarray}
\psi_1^\alpha(x,h_n) & = &\frac{1}{V_\alpha(x,h_n)^{1/4}}
\exp\left(-\frac{1}{h_n}\int_{x_+^\alpha(h_n)}^x\sqrt{V_\alpha(z,h_n)}~dz\right) \nonumber \\
 & & \times \left(1+\sum_{j=1}^{j_0}u_j^\alpha(x)h_n^j+R_{j_0+1}(x,h_n)\right), \label{WKBcubHn}
\end{eqnarray}
where $|u_j^\alpha(x)| = \mathcal{O}(|x|^{-5j/2})$ and $|R_{j_0+1}(x,h)|\leq C|x|^{-5(j_0+1)/2}h^{j_0+1}$.
\end{theorem}
In particular the expansion (\ref{WKBcubHn}) holds uniformly for $x\in\mathbb{R}$.

\begin{figure}[H]\begin{center}
 \scalebox{0.3}{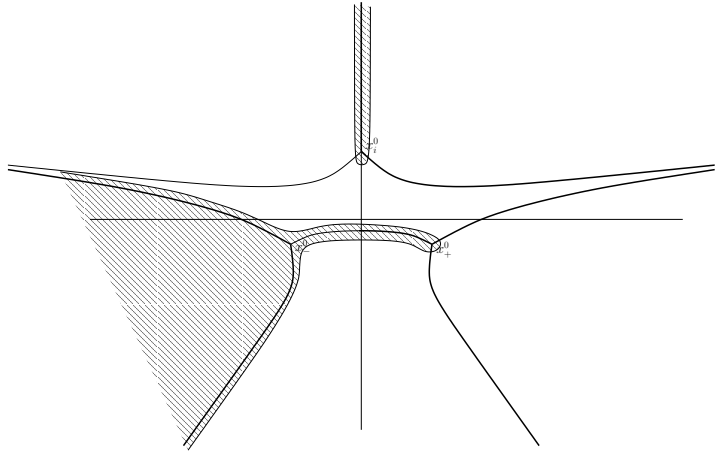}
 \caption{The domain $\Lambda_+(\varepsilon)$ (unshaded domain). $\Lambda_-(\varepsilon)$ is obtained from $\Lambda_+(\varepsilon)$ by applying 
 the symmetry
 of axis $i\mathbb{R}$.}
\label{figLambda}
\end{center}
\end{figure}

\textbf{Proof: } 
We apply Theorem $3.1$, ch. $10$, p. $366$ of \cite{Olv}.\\
Let
\[
 S(x) = \int_{x_+^0}^x\sqrt{iz^3-1}~dz,
\]
where $x_+^0=x_+^\alpha(0)$, and let $\Lambda_\pm$ be the set of points $x\in\mathbb{C}$ such that there exists a path
$\gamma_x$ joining $\pm\infty$ to $x$ such that $\Re S\circ\gamma_x$ is increasing (\emph{canonical path}).
Let $\Lambda_\pm(\varepsilon) = \{x\in\Lambda_\pm : d(x,\pa\Lambda_\pm)\geq\varepsilon\}$ (see Figure \ref{figLambda}).
We then notice that
\[
 \Gamma_\varepsilon = \Lambda_+(\varepsilon)\cup\Lambda_-(\varepsilon).
\]
According to Theorem $3.1$, ch. $10$, p. $366$ of \cite{Olv},
there exists $h_0>0$ such that, for $h\in]0,h_0[$, any solution $\psi_\pm^\alpha(\cdot,h)\in L^2(\mathbb{R}^\pm)$
satisfies (\ref{WKBcubXn}) and (\ref{WKBcubHn}) in $\Lambda_\pm(\varepsilon)$, up to a multiplicative constant
 $c_\pm(h)\in\mathbb{C}$, and with $h\rightarrow0$ instead of the sequence $(h_n)_n$. Indeed,
in order to check that the bound $(3.04)$ in \cite{Olv} on the remainder term of order $k$ is of size $\mathcal{O}(h^k)$,
we check that the conditions $(i)-(iv)$ p. $370$
are satisfied, which can be done by observing that the function
\[
\sigma_\alpha(x,h) := \frac{1}{V_\alpha(x,h)^{3/4}}\left[\frac{1}{V_\alpha(x,h)^{1/4}}\right]''
\]
satisfies, for some $k>0$,
\[
 |\sigma_\alpha(x,0)|\leq \frac{k}{1+|x|^5}~~~\textrm{ and }~~~\sigma_\alpha(x,h) = \sigma_\alpha(x,0)(1+\mathcal{O}(h^{4/5}))
\]
uniformly for $x\in\Lambda_\pm(\varepsilon)$.\\
To conclude, we have seen in Subsection \ref{pScaling} that if $\lambda_n(\alpha)$ denotes the $n$-th eigenvalue of
$\CA_\alpha$, and if
\beq\label{defHn}
h_n = \lambda_n(\alpha)^{-5/6},
\eeq
then there exists, for all $n\geq1$, a solution $\psi_1^\alpha(\cdot,h_n)\in L^2(\mathbb{R})$ of (\ref{solPsin}).
Then, according to the previous arguments, $\psi_1^\alpha(\cdot,h_n)$ satisfies (\ref{WKBcubXn}) and (\ref{WKBcubHn}) in
$\Lambda_+(\varepsilon)$ and $\Lambda_-(\varepsilon)$ up to respective constants $c_+(h)$ and $c_-(h)$.
Comparing these expressions for $x\in\Lambda_+(\varepsilon)\cap\Lambda_-(\varepsilon)$,
we see that $c_+(h) = c_-(h)$, and the statement follows by choosing $c_+(h) = c_-(h) = 1$.
\hfill $\square$\\

The asymptotic expansion (\ref{WKBcubHn}) does not hold in the neighborhood of the bounded Stokes line $\ell_f^\alpha(h)$.
In order to determine the behavior of a solution on $\ell_f^\alpha(h)$, we have to take into account the presence of terms
of the form
\[
 V_\alpha(x,h)^{-1/4}\exp\left(+\frac{1}{h}\int_{x_+^\alpha(h)}^x\sqrt{V_\alpha(z,h)}~dz\right)
\]
in its expression.
Those terms, exponentially small as \\$h^{-1}\Re\int_{x_+^\alpha(h)}^x\sqrt{V_\alpha(z,h)}~dz\rightarrow-\infty$,
are significant on $\ell_f^\alpha(h)$. In the following subsection,
we consider solutions which oscillate along $\ell_f^\alpha(h)$.
We will obtain an asymptotic expression which also holds in a neighborhood of the turning points $x_\pm^\alpha(h)$.

\subsection{Behavior of the eigenfunctions in the neighborhood of the turning points}\label{pPtsTour}
\begin{figure}[H]\begin{center}
 \scalebox{0.25}{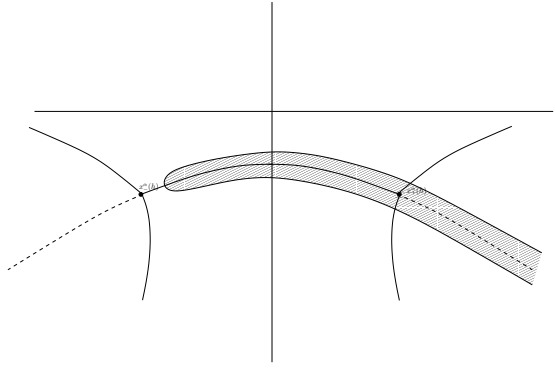}
 \caption{The domain $\CD_+(\delta,\eta)$ (shaded domain). The line joining $x_-^\alpha(h)$ to $x_+^\alpha(h)$ is the finite Stokes line
 $\ell_f^\alpha(h)$. The two dashed lines represent the anti-Stokes lines $\tilde\ell_-^\alpha(h)$ (on the left) and
 $\tilde\ell_+^\alpha(h)$ (on the right).}
\label{figCubicDom}
\end{center}
\end{figure}

In the neighborhood of a turning point, the previous asymptotic expansions are no longer available.
We will now use an approximation of the solutions involving the Airy function $Ai$.\\
We introduce the \emph{anti-Stokes lines} starting from $x_\pm^\alpha(h)$, defined as the level curves of the function
\[
 x\mapsto \Im \int_{x_+^\alpha(h)}^x\sqrt{V_\alpha(z,h)}~dz
\]
containing $x_\pm^\alpha(h)$.
A local analysis near the turning points shows that there exist three anti-Stokes lines starting from
$x_\pm^\alpha(h)$, and we will denote by $\tilde\ell_\pm^\alpha(h)$ (see Figure \ref{figCubicDom})
the one that satisfies
\[
 \forall x\in\tilde\ell_\pm^\alpha(h),~~ \int_{x_\pm^\alpha(h)}^x\sqrt{V_\alpha(z,h)}~dz>0.
\]
As in the previous subsection, we define a neighborhood of the line $\tilde\ell_\pm^0(0)$ by
\beq
\tilde\ell_{\pm,\delta}^0 = \{x\in\mathbb{C} : d(x,\tilde\ell_\pm^0(0))<\delta\},
\eeq
and we have $\tilde\ell_\pm^\alpha(h)\subset\tilde\ell_{\pm,\delta}^0$ for $h$ small enough.\\
Let $\eta>0$ be such that $\eta<|x_+^0(0)-x_-^0(0)|$.
Note that, for $h$ small enough, it implies $\eta<|x_+^\alpha(h)-x_-^\alpha(h)|$.
Then, for $\delta>0$, we denote
\beq
\CD_\pm(\delta,\eta) = \left(\ell_{f,\delta}^0\cap\{x\in\mathbb{C} : |x-x_\pm^0(0)|<\eta\}\right)\cup\tilde\ell_{\pm,\delta}^0.
\eeq
This domain is represented on Figure \ref{figCubicDom}.\\
In the following statement and its proof, we use the notation
\beq
 \zeta_\pm^\alpha(x,h) = \left(\frac{3}{2}\int_{x_\pm^\alpha(h)}^x\sqrt{V_\alpha(z,h)}dz\right)^{2/3}
\eeq
and
\beq\label{defPsiErr}
\tilde\sigma_\pm = 
\frac{1}{|V_\alpha|^{1/4}}\pa_x^2\left(\frac{1}{|V_\alpha|^{1/4}}\right)-\frac{5|V_\alpha|^{1/2}}{16|\zeta_\pm^\alpha|^3},
\eeq
which is defined for $x\neq x_\pm^\alpha(h)$.

\begin{theorem}\label{LemVoisZeros}
 Let $\alpha\in\mathbb{R}$. There exist positive constants $\delta>0$ and $h_1>0$, and two solutions $\psi_\pm^\alpha(x,h)$ of
equation
\[
 \CA_\alpha(h)\psi_\pm^\alpha(x,h) = \left(-h^2\frac{d^2}{dx^2}+V_\alpha(x,h)\right)\psi_\pm^\alpha(x,h) = 0
\]
such that, for all $h\in(0,h_1]$ and $x\in \CD_\pm(\delta,\eta)$,
\beq\label{exprPsipm}
\psi_\pm^\alpha(x,h) = \left(\frac{\zeta_\pm^\alpha(x,h)}{V_\alpha(x,h)}\right)^{1/4}
Ai\left(\frac{\zeta_\pm^\alpha(x,h)}{h^{2/3}}\right)+hr_\pm^\alpha(x,h),
\eeq
where the function $r_\pm^\alpha$ satisfies, for all $h\in[0,h_1]$,
\beq\label{borner}
\left\{
\begin{array}{lll}
\forall x\in \CD_\pm(\delta,\eta)\setminus\ell_f^\alpha(h), & |r_\pm^\alpha(x,h)|\leq 
C_\pm^\alpha(x)\left|Ai\left(\frac{\zeta_\pm^\alpha(x,h)}{h^{2/3}}\right)\right|, & \\
\forall x\in\CD_\pm(\delta,\eta)\cap\ell_f^\alpha(h), & |r_\pm^\alpha(x,h)|\leq K_\pm^\alpha, & 
\end{array}\right.
\eeq
for some constant $K_\pm^\alpha>0$ and some function
$C_\pm^\alpha(x)$ bounded in $\CD_\pm(\delta_\pm,\eta)$ outside any open neighborhood of
$\ell_f^0(0)\,$.
\end{theorem}
\textbf{Proof: } 
We work in the domain $\CD_+(\delta,\eta)$, and we will possibly drop the index
$+$ in the expressions.
We shall apply Theorem $9.1$, p. $417$ in \cite{Olv}, with a $h$-dependent potential here.
We introduce the following change of variable in $\CD_+(\delta,\eta)$ 
($\delta$ small enough will be determined in the following):
\beq\label{xtozeta}
 x\mapsto \zeta = \zeta(x,h)
\eeq
for a fixed $h\in[0,h_0]$. We denote its inverse by
\beq\label{zetatox}
 \zeta\mapsto x = x(\zeta,h).
\eeq
The three Stokes lines starting from $x_+^\alpha(h)$ are mapped by (\ref{xtozeta}) onto the half-lines
\[
 L_j=\arg^{-1}\left\{\frac{\pi}{3}+\frac{2j\pi}{3}\right\},
\]
and the anti-Stokes line $\tilde\ell_+^\alpha(h)$ is mapped onto the half-line $[0,+\infty[$.\\

Let $a=+\infty$, and let $Z(a)$ 
be the set of points $\zeta\in\mathbb{C}$ such that there exists a complex path $\gamma_\zeta$ joining $\zeta$ to $a$, 
which coincides at infinity with $[0,+\infty[$, and such that $v\mapsto\Re\gamma_\zeta(v)^{3/2}$ is
non-decreasing.\\
Then there exists $\delta>0$ such that, for $h=0$, $\zeta(\CD_+(2\delta,\eta),0)\subset Z(a)$.
Since $V_\alpha$ has the form 
\beq\label{formValpha}
V_\alpha(x,h) = V_0(x)+h^{4/5}v_\alpha(x,h),
\eeq
where $|v_\alpha(x,h)| = o(|V_0(x)|)$
uniformly with respect to $h$ as $|x|\rightarrow+\infty$, there exists $h_1>0$ such that for all $h\in]0,h_1[$,
\[
 \zeta(\CD_+(\delta,\eta),h)\subset Z(a).
\]
Thus, Theorem $9.1$, p. $417$ in \cite{Olv}, which applies for all $\zeta\in Z(a)$, ensures that there exists a solution
\[
 \psi^\alpha(x,h) = \left(\frac{\zeta(x,h)}{V_\alpha(x,h)}\right)^{1/4}W(\zeta(x,h),h),
\]
where $W$ has the form
\beq\label{solW}
\forall h\in(0,h_1],~\forall \zeta\in\zeta(\CD(\delta,\eta),h),~~
W(\zeta,h) = Ai\left(\frac{\zeta}{h^{2/3}}\right)+h\varepsilon(\zeta,h).
\eeq
In view of inequality ($9.03$), p. $418$ in \cite{Olv} (here applied with $n=0$, $u=h^{-1}$ and $\varepsilon_{2n+1}$ replaced by $h\varepsilon(\zeta,h)$),
in order to prove that the function
\beq\label{defR}
 r_+^\alpha(x,h) :=\left(\frac{\zeta(x,h)}{V(x,h)}\right)^{1/4}\varepsilon(\zeta(x,h),h)
\eeq
satisfies the bounds (\ref{borner}), it remains to check that there exists $M>0$ such that, for all $h\in]0,h_1[$ and
$\zeta\in\zeta(\CD_+(\delta,\eta),h)$,
\beq\label{intErr}
 \int_{x(\gamma_\zeta,h)}|\tilde\sigma(z,h)||dz| \leq M,
\eeq
where $\tilde\sigma$ is the function defined in (\ref{defPsiErr}), and $x(\gamma_\zeta,h)$ denotes the image by (\ref{zetatox})
of the path $\gamma_\zeta$ defined above.
Here we used the notation $|dz| = |x(\gamma_\zeta,h)'(t)|dt$.\\
Notice that the function $\tilde\sigma(x,h)$ is integrable at $x=x_\pm^\alpha(h)$, see for instance Lemma $3.1$, p. $399$ in \cite{Olv}.
Moreover, one can easily check that there exists $k>0$ such that
\beq\label{borneE}
 |\tilde\sigma(x,0)|\leq\frac{k}{1+|x|^{7/2}}
\eeq
for $|x|$ large enough, $x\in\CD_+(\delta,\eta)$. Thus, (\ref{intErr}) follows from (\ref{formValpha}) and (\ref{borneE}), and (\ref{borner})
is then proved.
\hfill $\square$\\

We now want to integrate the solution $\psi_\pm^\alpha$ over a path on which $\zeta(x,h)$ is real. 
In this purpose, we choose a $\mathcal{C}^1$ path
$\gamma_h = \gamma_{h,\pm}^\alpha : [-d,+\infty[\rightarrow\mathbb{C}$ such that $\gamma_h(0) = x_\pm^\alpha(h)$,
\beq\label{imGamma}
 \gamma_h([-d,+\infty[) = \bar\CD_\pm(\delta,\eta)\cap(\ell_f^\alpha(h)\cup\tilde\ell_\pm^\alpha(h)),
\eeq
and satisfying
\beq\label{gammaReg}
\forall t\in[-d,+\infty[,~~|\gamma_h'(t)| = 1.
\eeq
Such a smooth path exists because both lines $\ell_f^\alpha(h)$ and $\tilde\ell_\pm^\alpha(h)$ reach 
the point $x_\pm^\alpha(h)$ with the same angle
$-\frac{2}{3}\arg\sqrt{\pa_xV_\alpha(x_\pm^\alpha(h),h)}$ (modulo $\pi$).\\
Let us fix $\delta'\in]0,\delta[$, $\eta'\in]0,\eta[$, and
$\chi_\pm\in\mathcal{C}^\infty(\mathbb{C},[0,1])$ with $\chi_\pm(x)=1$ for $x\in\CD_\pm(\delta',\eta')$ and
$\Supp\chi_\pm\subset\CD_\pm(\delta,\eta)$.
\begin{lemma}\label{Corw}
There exists $c_\pm^\alpha\neq0$ such that, as $h\rightarrow0$,
\beq\label{normCubZeros}
\int_{\gamma_{h,\pm}^\alpha}\psi_\pm^\alpha(x,h)^2\chi_\pm^\alpha(x)dx = c_\pm^\alpha h^{1/3}(1+o(1)).
\eeq
\end{lemma}
\textbf{Proof: } 
Let us consider the case of $\psi_+^\alpha$. We set $\zeta=\zeta_+^\alpha$, $\gamma_h = \gamma_{h,+}^\alpha$,
$\chi=\chi_+^\alpha$ to simplify the notation.\\
We first apply the following change of variable, for a fixed $h\in[0,h_1]$:
\[
 [-d,+\infty[\ni t\mapsto \zeta := \zeta(\gamma_h(t),h)\in[-b_h,+\infty[,
\]
where $[-b_h,+\infty[$ is the range of this function.
Note that we have $\gamma_h(t) = x(\zeta,h)$, where $x(\cdot,h)$ is the inverse mapping (\ref{zetatox}).\\
Let $b$ such that $b>b_h$ for all $h\in[0,h_1]$, and $\chi_h(\zeta) = \chi\circ x(\zeta,h)$, supported in
$]-b,+\infty[$. Then,
\beq\label{decompNormW}
\int_{\gamma_h}\psi_+^\alpha(x,h)^2\chi(x)dx = I_0(h)+hI_1(h)+h^2I_2(h)
\eeq
where
\beq\label{I0w}
 I_0(h) = \int_{-b}^{+\infty}\frac{\zeta}{V_\alpha(x(\zeta,h),h)}Ai\left(\frac{\zeta}{h^{2/3}}\right)^2\chi_h(\zeta)d\zeta,
\eeq
\beq\label{I1w}
I_1(h) = 2\int_{-b}^{+\infty}
\frac{\zeta}{V_\alpha(x(\zeta,h),h)}Ai\left(\frac{\zeta}{h^{2/3}}\right)\varepsilon(\zeta,h)\chi_h(\zeta)d\zeta
\eeq
and
\beq\label{I2w}
I_2(h) = \int_{-b}^{+\infty}\frac{\zeta}{V_\alpha(x(\zeta,h),h)}\varepsilon(\zeta,h)^2\chi_h(\zeta)d\zeta.
\eeq
We recall that the Airy function is defined by
\[
 Ai(x) = \frac{1}{2\pi}\int_\mathbb{R}e^{i(x\xi+\xi^3/3)}d\xi,
\]
hence
\[
 Ai\left(\frac{\zeta}{h^{2/3}}\right) = \frac{1}{2\pi h^{1/3}}\int_\mathbb{R}e^{\frac{i}{h}(\zeta\xi+\xi^3/3)}d\xi.
\]
Thus,
\beq\label{I0iii}
 I_0(h) = \frac{1}{4\pi^2h^{2/3}}\iiint_{[-d,+\infty[\times\mathbb{R}^2}\frac{\zeta}{V_\alpha(x(\zeta,h),h)}
e^{\frac{i}{h}\Phi(\zeta,\eta,\xi)}\chi_h(\zeta)d\zeta d\eta d\xi,
\eeq
where
\[
 \Phi(\zeta,\eta,\xi) = \zeta(\xi-\eta)+\frac{1}{3}(\xi^3-\eta^3).
\]
It is then straightforward to check that for all $\xi\in\mathbb{R}$,
the function $\Phi(\cdot,\cdot,\xi)$  has a unique critical point $(-\xi^2,\xi)$, which is non-degenerate.
Moreover, $\Phi(-\xi^2,\xi,\xi) = 0$. 
Thus, the stationary phase method with $\xi$ fixed in (\ref{I0iii}), yields
\beq\label{I0wfinal}
I_0(h) = c_+^\alpha h^{1/3}(1+o(1)),~~h\rightarrow0,
\eeq
where
\beq\label{exprC0}
 c_+^\alpha = -(2\pi)^{-3/2}\int_{\xi\in\mathbb{R}}\frac{\xi^2}{V_\alpha(x(-\xi^2,0),0)}\chi_0(-\xi^2)d\xi.
\eeq
Finally, using (\ref{borner}) and the asymptotic behavior of the Airy function as $z\rightarrow\pm\infty$ (see
\cite{AbrSteg}), one can easily check that
\[
 hI_1(h)+h^2I_2(h) = \mathcal{O}(h^{7/6}),
\]
and the statement follows.
\hfill $\square$\\

\subsection{Connection}\label{pRaccord}
In Subsections \ref{pWKB} and \ref{pPtsTour}, we have determined the asymptotic behavior as $h\rightarrow0$ of
several solutions of (\ref{probCubic}). More precisely, we have built a solution
$\psi_1^\alpha(\cdot,h_n)\in L^2(\mathbb{R})$
whose behavior is known in a domain
$\Gamma_\varepsilon$ avoiding a neighborhood of the bounded Stokes line $\ell_f^\alpha(h)$,
and two solutions $\psi_\pm^\alpha(\cdot,h)$ whose asymptotic behavior is known in a neighborhood of
$\ell_f^\alpha(h)$ avoiding the opposite turning point (see Theorem \ref{LemVoisZeros}).
We now want to connect these solutions, comparing their asymptotic expressions in the intersection of their domain of
validity.\\
We first state the  Bohr-Sommerfeld quantization rule, which gives a relation between the value of
$h_n$ and the index $n$. We will then use it to determine the coefficient relating the solutions
$\psi_1^\alpha$ and $\psi_\pm^\alpha$. This lemma can be proved as Formula $(25)$ in \cite{GreMaiMar09}.
\begin{lemma}[Bohr-Sommerfeld quantization rule]
\beq\label{exprBS}
\Im\int_{x_-^\alpha(h_n)}^{x_+^\alpha(h_n)}\sqrt{V_\alpha(z,h_n)}~dz = \pi\left(n+\frac{1}{2}\right)h_n+\mathcal{O}(h_n^2).
\eeq
\end{lemma}
We are now going to compare the asymptotic expressions of $\psi_1^\alpha$ and $\psi_\pm^\alpha$, for fixed $h$
as $|x|\rightarrow+\infty$
along the lines $\tilde\ell_\pm^\alpha(h)$. Let $n\geq1$ be large enough so that
$\tilde\ell_-^\alpha(h_n)\subset\tilde\ell_{-,\delta}^0$, and let $x\in\tilde\ell_-^\alpha(h_n)$.
We are then able to use the asymptotic expansion of the Airy function as $|z|\rightarrow+\infty$ \cite{AbrSteg}, $|\arg z|<\pi$,
with $z = \zeta_-^\alpha(x,h_n)$. If we denote $S_\pm^\alpha(x,h) = \int_{x_\pm^\alpha(h)}^x\sqrt{V_\alpha(z,h)}~dz$,
expression (\ref{exprPsipm}) then writes
\begin{eqnarray}
 \psi_-^\alpha(x,h_n) & = & \frac{h_n^{1/6}}{2\sqrt\pi V_\alpha(x,h_n)^{1/4}}\exp\left(-\frac{1}{h_n}S_-^\alpha(x,h_n)\right)
 (1+\mathcal{O}(S_-^\alpha(x,h_n)^{-3/2})) \nonumber\\
  & = & \frac{h_n^{1/6}}{2\sqrt\pi}\exp\left(-\frac{1}{h_n}\int_{x_-^\alpha(h_n)}^{x_+^\alpha(h_n)}\sqrt{V_\alpha(z,h_n)}~dz\right)\psi_1^\alpha(x,h_n)
  \nonumber\\
   & & \times(1+\mathcal{O}(|x|^{-5/2})), \label{lien1-}
\end{eqnarray}
where we used (\ref{WKBcubHn}).\\

The two solutions $\psi_-^\alpha$ and $\psi_1^\alpha$ being both exponentially decreasing as
$|x|\rightarrow+\infty$ along $\tilde\ell_-^\alpha(h_n)$, they are necessarily colinear.
Hence, (\ref{exprBS}) and (\ref{lien1-}) yield
\beq\label{racc1-}
\psi_-^\alpha(x,h_n) = \frac{(-1)^{n-1}i}{2\sqrt\pi}h_n^{1/6}\psi_1^\alpha(x,h_n)(1+\mathcal{O}(h_n)),~~~n\rightarrow+\infty.
\eeq
Similarly, comparing the asymptotic representations of $\psi_1^\alpha$ and $\psi_+^\alpha$ as $|x|\rightarrow+\infty$ along
$\tilde\ell_+^\alpha(h_n)$, we get
\beq\label{racc1+}
\psi_+^\alpha(x,h_n) = \frac{1}{2\sqrt\pi}h_n^{1/6}\psi_1^\alpha(x,h_n).
\eeq

Due to these relations, we can integrate the square of the solution
$\psi_1^\alpha(x,h_n)$ over the curve consisting in the union of the three lines
$\tilde\ell_-^\alpha(h_n)$, $\ell_f^\alpha(h_n)$ and $\tilde\ell_+^\alpha(h_n)$,
\beq\label{defLalpha}
 \mathcal{L}_\alpha(h_n) = \tilde\ell_-^\alpha(h_n)\cup\ell_f^\alpha(h_n)\cup\tilde\ell_+^\alpha(h_n).
\eeq
We choose $\eta>0$ such that $\eta<|x_+^0(0)-x_-^0(0)|$ and such that
$\ell_{f,\delta}^0\subset\CD_+(\delta,\eta)\cup\CD_-(\delta,\eta)$.
Let also $\eta'<|x_+^0(0)-x_-^0(0)|/2$ and $\delta'\in]0,\delta[$.
We choose a partition of unity $(\chi_-,\chi_+)$ such that,
for all $h\in]0,h_1]$ and all $x\in\mathcal{L}_\alpha(h)$, $\chi_-(x)+\chi_+(x) = 1$, and such that
$\chi_\pm(x) = 1$ for $x\in\CD_\pm(\delta',\eta')$, and $\Supp\chi_\pm\subset\CD_\pm(\delta,\eta)$.\\
Then, according to (\ref{racc1+}) and
(\ref{racc1-}), for all $x\in\mathcal{L}_\alpha(h_n)$,
\beq\label{decompPsi1}
\psi_1^\alpha(x,h_n)^2 = 4\pi h_n^{-1/3}(\psi_+^\alpha(x,h_n)^2\chi_+(x)-\psi_-^\alpha(x,h_n)^2\chi_-(x))
 (1+\mathcal{O}(h_n))
\eeq
as $n\rightarrow+\infty$.\\
Thus, we deduce the following lemma from (\ref{normCubZeros}), where $c_\alpha = c_+^\alpha+c_-^\alpha\neq0$
(see (\ref{exprC0})):
\begin{lemma}\label{lemIntPsi1}
For all $\alpha\in\mathbb{R}$, there exists $c_\alpha\neq0$ such that
\beq\label{normL}
\int_{\mathcal{L}_\alpha(h_n)}\psi_1^\alpha(x,h_n)^2dx = c_\alpha(1+o(1))
\eeq
as $n\rightarrow+\infty$.
\end{lemma}

In the last section, we gather the previous results to prove Theorem \ref{mainCubic}.

\section{Estimate on the instability indices}\label{sCcl}
Let us first recall the following general result, which will provide an explicit formula for the instability indices
$\kappa_n(\CA_\alpha)$, for $n$ large enough (see \cite{AslDav}).

\begin{lemma}\label{propII}
Let $\CA$ be a closed operator on the Hilbert space $\CH$, and $\lambda\in\sigma(\CA)$ a simple isolated eigenvalue.
Let $\Pi_\lambda$ be the spectral projectioon associated with $\lambda$, $u_\lambda$ an eigenvector associated with $\lambda$,
and $u_\lambda^*$ an eigenvector of $\CA^*$ associated with the eigenvalue $\bar\lambda$. Then:
\begin{enumerate}[(i)]
 \item $\Pi_\lambda$ has rank $1$ if and only if $\sca{u_\lambda}{u_\lambda^*}\neq0$.
\item In this case, we have
\beq\label{exprII}
\kappa(\lambda):=\|\Pi_\lambda\| = \frac{\|u_\lambda\|\|u_\lambda^*\|}{|\sca{u_\lambda}{u_\lambda^*}|}.
\eeq
\end{enumerate}
\end{lemma}

We recall (see Subsection \ref{pScaling}) that the eigenfunctions $u_n^\alpha$
associated with the $n$-th eigenvalue $\lambda_n(\alpha)\in\mathbb{R}$ of $\CA_\alpha$ have the form
\beq
 u_n^\alpha(x) = \psi_\alpha(h_n^{2/5}x,h_n),
\eeq
where
\beq\label{RappelHn}
h_n = \lambda_n(\alpha)^{-5/6},
\eeq
and where $\psi_\alpha(\cdot,h_n)\in L^2(\mathbb{R})$ is a solution of $\CA_\alpha(h_n)\psi_\alpha(\cdot,h_n) = 0$.\\
We normalize $u_n^\alpha$ so that
\beq\label{RappelUn}
 u_n^\alpha(x) = \psi_1^\alpha(h_n^{2/5}x,h_n),
\eeq
where $\psi_1^\alpha$ is the solution introduced in Theorem \ref{CorWKB}.\\
We have
\begin{proposition}\label{propDenomCubic}
Let $\alpha\geq0$. There exists $N\geq1$ such that, for all $n\geq N$, the spectral projection $\Pi_n(\alpha)$ of $\CA_\alpha$
associated with $\lambda_n(\alpha)$ has rank $1$. Moreover, there exists $k_\alpha>0$ such that the $n$-th
instability index satisfies
\beq\label{exprIIcubic}
\kappa_n(\alpha) = k_\alpha\|\psi_1^\alpha(\cdot,h_n)\|_{L^2(\mathbb{R})}(1+o(1)),~~~n\rightarrow+\infty.
\eeq
\end{proposition}
\textbf{Proof: } 
By deformation of the integration path, and using the exponential decay of
$\psi_1^\alpha(x,h_n)$ as $|x|\rightarrow+\infty$ in the sectors $\arg^{-1}(]-3\pi/10,\pi/10[)$ and
$\arg^{-1}(]9\pi/10,13\pi/10[)$
(see Theorem \ref{CorWKB}), we get
\beq\label{homotopy}
\int_\mathbb{R}\psi_1^\alpha(x,h_n)^2dx = \int_{\mathcal{L}_\alpha(h_n)}\psi_1^\alpha(x,h_n)^2dx.
\eeq
We then notice that $\CA_\alpha^*\Gamma = \Gamma\CA_\alpha$, where
$\Gamma : u(x)\mapsto\barr{u(x)}$. Hence, we have
$
(u_n^\alpha)^*(x) = \barr{u_n^\alpha(x)} 
$,
with the notation of Proposition \ref{propII}. Thus, according to (\ref{RappelUn}),
\[
\sca{u_n^\alpha}{(u_n^\alpha)^*}  =  h_n^{-2/5}\int_{\mathcal{L}_\alpha(h_n)}\psi_1^\alpha(x,h_n)^2dx.
\]
Using (\ref{normL}) we then get, for $n$ large enough,
$|\sca{u_n^\alpha}{(u_n^\alpha)^*}|>0$,
and the desired statement on the rank of $\Pi_n(\alpha)$ follows from Proposition \ref{propII}, $(i)$.
Expression (\ref{exprIIcubic}) follows from (\ref{normL}) and Proposition \ref{propII}, $(ii)$, after the change
of variable $x\mapsto h_n^{2/5}x$.
\hfill $\square$\\

Now it remains to determine an equivalent for the norm $\|\psi_1^\alpha(\cdot,h_n)\|_{L^2(\mathbb{R})}$ appearing in
(\ref{exprIIcubic}). We will do so by using the expansion (\ref{WKBcubHn}).
Let us recall that this expansion is uniform with respect to $x\in\mathbb{R}$, hence by integrating:
\beq\label{normPsi1}
\|\psi_1^\alpha(\cdot,h_n)\|_{L^2(\mathbb{R})}^2 = (1+o(1))\int_\mathbb{R}a(x)e^{-\varphi_\alpha(x,h_n)}dx
\eeq
as $n\rightarrow+\infty$,
where
\[
 a(x) = \frac{1}{V_0(x)^{1/4}}
\]
and
\[
 \varphi_\alpha(x,h) = \frac{2}{h}\Re\int_{x_+^\alpha(h)}^x\sqrt{V_\alpha(z,h)}~dz.
\]
\begin{lemma}\label{propNumerCubic}
If $\alpha\geq0$ then, as $n\rightarrow+\infty$,
\beq\label{eqNumerCubic}
 \|\psi_1^\alpha(\cdot,h_n)\|_{L^2(\mathbb{R})}^2 = \frac{\sqrt{2}}{2}\Gamma(1/4)h_n^{1/4}(1+o(1))
\exp\left(\frac{C}{h_n}+\frac{\alpha r}{h_n^{1/5}}\right),
\eeq
where
\beq\label{CetR}
C = \int_0^1\sqrt{1-t^3}~dt~>0~~\textrm{ and }~~ r = \frac{1}{2}\int_0^1\frac{t}{\sqrt{1-t^3}}~dt.
\eeq
\end{lemma}
\textbf{Proof: } 
Let us first assume that $\alpha>0$.
We shall apply the Laplace method with two parameters in \cite{Ped} to determine the behavior as $h\rightarrow0$
of the integral
\[
 I_\alpha(h) = \int_\mathbb{R}a(x)e^{-\varphi_\alpha(x,h)}dx
\]
appearing in (\ref{normPsi1}). We write
$\varphi_\alpha(x,h) = \frac{1}{h}g_\alpha(x,\varepsilon(h))$ with $\varepsilon(h) = h^{4/5}$ and
\[
 g_\alpha(x,\varepsilon) = 2\Re\int_{\tilde x_+^\alpha(\varepsilon)}^x\sqrt{\tilde V_\alpha(z,\varepsilon)}~dz,
\]
where we have denoted $\tilde x_+^\alpha(\varepsilon) = x_+^\alpha(\varepsilon^{5/4})$ and
$\tilde V_\alpha(x,\varepsilon) = V_\alpha(x,\varepsilon^{5/4})=ix^3+i\alpha\varepsilon x-1$.
\\
The function $g_\alpha$ is $\mathcal{C}^\infty$ for $x\in\mathbb{R}$ and $\varepsilon$ small enough.
Moreover, $g_\alpha(\cdot,0)$ has a unique critical point $x=0$. Indeed,
\[
 \pa_xg_\alpha(x,0) = 2\Re\sqrt{ix^3-1} = 0
\]
if and only if $\arg(ix^3-1) = \pi$, that is $x=0$.\\
We write
\beq
 \varphi_\alpha(x,h) = \frac{1}{h}g_\alpha(x,0)+\frac{\varepsilon(h)}{h}\pa_\varepsilon g_\alpha(x,0)+
\mathcal{O}\left(\frac{\varepsilon(h)^2}{h}\right),
\eeq
and we easily check that the remainder term is uniform with respect to $x\in\mathbb{R}$. We also check that
\[
 \pa_x^2g_\alpha(0,0) = \pa_x^3g_\alpha(0,0) = 0~~~\textrm{ and }~~~\pa_x^4g_\alpha(0,0) = 6,
\]
and that
\[
 \pa_x\pa_\varepsilon g_\alpha(0,0) = 0~~~\textrm{ and }~~~ \pa_x^2\pa_\varepsilon g_\alpha(0,0) = \alpha.
\]
Thus,
\[
 g_\alpha(x,0)-g_\alpha(0,0) = \frac{x^4}{4}+\mathcal{O}(|x|^5),~~
\pa_\varepsilon g_\alpha(x,0)-\pa_\varepsilon g_\alpha(0,0) = \frac{\alpha x^2}{2}
+\mathcal{O}(|x|^3).
\]
We can then apply Theorem $2$ in \cite{Ped}, with $\phi(x) = g_\alpha(x,0)-g_\alpha(0,0)$, 
$\psi(x) = -\pa_\varepsilon g_\alpha(x,0)$, $\nu = 4$, $\mu=2$, $\lambda=0$, and replacing $h$ by $h^{-1}$ and $k$
by $\varepsilon(h)h^{-1} = h^{-1/5}$. This yields
\[
 \|\psi_1^\alpha(\cdot,h_n)\|_{L^2(\mathbb{R})}^2 = \frac{\sqrt{2}}{2}\Gamma(1/4)h_n^{1/4}(1+o(1))
\exp\left(-\frac{1}{h_n}(g_\alpha(0,0)+h_n^{4/5}\pa_\varepsilon g_\alpha(0,0)\right)
\]
In order to get the desired statement, it only remains to notice that
$g_\alpha(0,0) = -C$ and $\pa_\varepsilon g_\alpha(0,0) = -\alpha r$, where $C$ and $r$ are the constants in (\ref{CetR}).\\
In the case $\alpha = 0$, we check similarly that the Laplace method applies (see for instance \cite{Erd})
and leads to the same statement.
\hfill $\square$\\

To conclude the proof of Theorem \ref{mainCubic},
we use the Bohr-Sommerfeld rule (\ref{exprBS}), which gives 
an asymptotic expansion for $h_n$. Let us compute the first few terms.
By expanding the left-hand-side of (\ref{exprBS}), we get
\[
 \Im\int_{x_-^\alpha(h_n)}^{x_+^\alpha(h_n)}\sqrt{V_\alpha(z,h_n)}~dz = \sqrt{3}C-\sqrt{3}\alpha rh_n^{4/5}+\mathcal{O}(h_n^{8/5}),
\]
where $C$ and $r$ are the constants in (\ref{CetR}). Expression (\ref{exprBS}) then writes
 \beq\label{exprBSbis}
 h_n = \frac{\sqrt{3}C}{\pi\left(n+\frac{1}{2}\right)}-\frac{3^{9/10}\alpha rC^{4/5}}{\pi^{9/5}\left(n+\frac{1}{2}\right)^{9/5}}
 +\mathcal{O}((n+1/2)^{-13/5}).
 \eeq
Gathering (\ref{exprIIcubic}), (\ref{eqNumerCubic}) and (\ref{exprBSbis}), and replacing $C$ and $r$ by their values
\[
 C = \frac{2\sqrt{3}\pi^{3/2}}{15\Gamma(2/3)\Gamma(5/6)}~~\textrm{ and }~~
r = \frac{\Gamma(2/3)\Gamma(5/6)}{2\sqrt{\pi}},
\]
we get the following statement, and Theorem \ref{mainCubic} follows.
\begin{theorem}
For all $\alpha\geq0$, there exists a positive constant $K_\alpha$ such that
\beq
\kappa_n(\alpha) = \frac{K_\alpha}{n^{1/4}}(1+o(1))\exp\left(\frac{\pi}{\sqrt{3}}n+\alpha c n^{1/5}\right)
\eeq
as $n\rightarrow+\infty$, where
\[
 c = (5/2)^{1/5}\pi^{-3/5}\Gamma(2/3)^{6/5}
\Gamma(5/6)^{6/5}.
\]
\end{theorem}

\textbf{Acknowledgments}\\
I am grateful to Bernard Helffer and André Martinez for their valuable help and comments. I acknowledge the support of the
ANR NOSEVOL.

\end{document}